\documentclass[reqno,oneside]{amsart}
\usepackage{times,amssymb,amsthm,eucal}
\usepackage{graphicx}
\newcommand{\half}{\tfrac{1}{2}}
\newcommand{\third}{\tfrac{1}{3}}
\newcommand{\egypt}{\tfrac{2}{3}}

\newcommand{\C}{\mathbb{C}}
\newcommand{\Z}{\mathbb{Z}}
\newcommand{\Q}{\mathbb{Q}}
\newcommand{\R}{\mathbb{R}}
\newcommand{\PS}{\mathbb{P}}
\newcommand{\OB}{\mathcal{O}}

\DeclareMathOperator{\rank}{rank}
\DeclareMathOperator{\SL}{SL}
\DeclareMathOperator{\Aut}{Aut}
\DeclareMathOperator{\Bir}{Bir}

\newcommand{\Aone}{{\mathcal{A}\mathrm{mp}}}

\newcommand{\Aclone}{{\,\overline{\!\mathcal{A}\mathrm{mp}}}}
\newcommand{\Mone}{{\mathcal{M}\mathrm{ov}}}
\newcommand{\Mcone}{{\mathcal{M}\mathrm{ov}\!_\smallplus}}
\newcommand{\Mclone}{{\overline{\mathcal{M}\mathrm{ov}}}}
\newcommand{\Weilsp}{{\widetilde N^1}}
\newcommand{\Cartsp}{{N^1}}

\newcommand{\HMB}{\mathcal{F}}
\newcommand{\HGp}{{H_5}}
\newcommand{\muGp}{{\mu_5}}

\newcommand{\Cox}{\mbox{$\circ\!\!=\!\!\circ\!\!=\!\!\circ\!\!=\!\!\circ$}}
\newcommand{\MCox}{\mbox{$\bullet\!\!=\!\!\circ\!\!=
  \!\!\circ\!\!=\!\!\bullet$}}
\newcommand{\TilCox}{\mbox{$\bullet\!\!=\!\!\circ\!\!=\!\!\circ$}}
\newcommand{\CubCox}{\mbox{$\bullet\!\!=\!\!\circ\quad\bullet$}}

\newcommand*{\stackrels}[3]{\mathrel{\mathop{#3}\limits^{#1}_{#2}}}
\newcommand{\ratto}{\mathrel{\relbar\!\:\to}}
\newcommand{\ratgets}{\mathrel{\gets\!\:\relbar}}
\newcommand{\rattogets}{\mathrel
  {\vcenter{\baselineskip0pt\lineskip0pt\hbox{$\ratto$}\hbox{$\ratgets$}}}}
\newcommand{\defeq}{\mathrel{\raisebox{.1ex}{:}}=}

\newcommand{\transpose}{{\scriptscriptstyle T}}
\newcommand{\smallplus}{{\scriptscriptstyle +}}
\newcommand{\dual}{\vee}

\newcommand{\ldot}{\mathbin{.}}

\let\secsign\S
\newcommand{\CY}{Cal\-abi\textup{--}\hspace{-.15em}Yau}
\newcommand{\HM}{Hor\-rocks\textup{--}Mum\-ford}

\newtheorem{theorem}{Theorem}
\newtheorem{lemma}[theorem]{Lemma}
\newtheorem{corollary}[theorem]{Corollary}
\theoremstyle{definition}
\newtheorem*{definition}{Definition}
\newtheorem*{convention}{Convention}

\begin{document}
\title{The Movable Fan of the Horrocks--Mumford Quintic}
\author{Michael Fryers\\Institut f\"ur Mathematik\\Universit\"at Hannover}
\thanks{Partially supported by
the Israel Academy of Sciences and Humanities
Centres of Excellence Programme,
the Emmy Noether Institute for Mathematics,
and the Minerva Foundation of Germany.}
\begin{abstract}
This paper explores the birational geometry
of a general \HM\ quintic threefold,
describing the set of all minimal models
up to marked isomorphism,
the movable fan (the way in which the
nef cones of all these models are arranged to form the movable cone),
and the birational automorphism group.
There are infinitely many models up to marked isomorphism,
falling into just eight (unmarked-) isomorphism classes.
\end{abstract}
\maketitle

The \HM\ bundle is a rank-$2$ vector bundle over $\PS^4$,
characterized by the fact that a general section vanishes along an Abelian
surface.  A general $2$-dimensional subspace of its sections defines a
quintic threefold with $100$ nodes.
Such threefolds (called \emph{\HM\ quintics}) have been studied
by Aure \cite{Aure89} and by Borcea \cite{Borcea91}.

Borcea described the $14$ projective small resolutions
of a general \HM\ quintic,
which are all smooth \CY\ threefolds, and also described an involution of
the moduli space of general \HM\ quintics which gives, for each such quintic,
a birationally equivalent quintic.
Up to isomorphism, there are just eight different threefolds
among the resolutions of these two quintics, which lie in
three different deformation families.

Here we describe the set of all minimal models of
a general \HM\ quintic up to \emph{marked} isomorphism,
the \emph{movable fan}, and the action of the 
birational automorphism group of the quintic on this fan.
It turns out that every one of these minimal
models is isomorphic to one of the eight described by
Borcea, although there are infinitely many up to marked isomorphism.
We check that Morrison's Movable Cone Conjecture
\cite{Morrison66} holds for these varieties.
The movable cone has a surprising form; 
its boundary is once but not twice differentiable.
In Section \ref{sec:detquintics} we describe the simpler movable cone of
determinantal quintics, of which \HM\ quintics are a specialisation.

This paper gives an example of a method which translates the
algebraic geometry of a variety into the real geometry and combinatorics of
the cones, and then continues the analysis purely in those terms, finally
recovering algebraic-geometric information from the description obtained.

\subsection*{Terminology: the movable fan}

A \emph{convex cone} in a real vector space is a convex subset closed
under multiplication by positive real scalars.

For a normal complex projective variety $X$,
we say that two divisors $D_1$ and $D_2$ on $X$ are
\emph{$\Q$-algebraically equivalent} if some non-zero multiple of
$(D_1-D_2)$ is algebraically equivalent to $0$.
For Cartier divisors, this is the same as numerical equivalence.
We define
\begin{equation*}
\Weilsp(X)=\{\,\mbox{divisors on $X$, modulo $\Q$-algebraic equivalence}\,\}
\otimes\R,
\end{equation*}
a finite-dimensional real vector space, and write $\Cartsp(X)$ for the
subspace generated by Cartier divisors.
Intersection with a curve in $X$ gives a linear form on $\Cartsp(X)$, but
this does not extend naturally to a linear form on $\Weilsp(X)$.

The \emph{nef cone} of $X$, written $\Aclone(X)$, is the closed
convex cone in $\Cartsp(X)$ consisting of all Cartier $\R$-divisors that
have non-negative intersection with every curve in $X$.
It spans the space $\Cartsp(X)$, and its relative interior is the
\emph{ample cone}, $\Aone(X)$, the convex cone generated by ample divisors.

A divisor on $X$ is said to be \emph{movable} if it is effective and
moves in a linear system whose base locus is of codimension at least $2$.
The \emph{movable cone} of $X$, written $\Mcone(X)$, is the convex hull
of the set of movable divisors on $X$; it is also a convex cone,
which spans $\Weilsp(X)$, and its interior, written $\Mone(X)$,
contains $\Aone(X)$.

If $Y$ is another normal projective variety, a birational equivalence
between $X$ and $Y$ is called an \emph{isomorphism in codimension $1$}
if it gives an isomorphism between open subsets of $X$ and $Y$ whose
complements are of codimension at least $2$.
Such a map gives an isomorphism of the divisor spaces $\Weilsp(X)$
and $\Weilsp(Y)$.
Under this isomorphism, $\Mcone(X)$ is carried to $\Mcone(Y)$,
but $\Aone(X)$ and $\Aone(Y)$ coincide only if the birational
equivalence is a biregular isomorphism, and are otherwise disjoint.
For \CY\ threefolds, it is known (from the existence and
termination of flops: see \cite{Kawamata88})
that the ample cones of all codimension $1$ isomorphic models of $X$
form a locally finite disjoint cover of $\Mone(X)$.
This arrangement of ample cones of various dimensions
covering the single movable
cone is called the \emph{movable fan}.

\section{\HM\ quintics}\label{sec:hmquintics}

In this section we review the definitions and basic properties
of the \HM\ quintics.
See \cite{HorMum73}, \cite{Borcea91}, and \cite{Aure89} for more details
of these constructions.

We recall \cite{HorMum73} that the \emph{\HM\ bundle},
$\HMB$, is a rank-$2$ vector bundle on $\PS^4$, of degree $5$.
A generic section of $\HMB$ vanishes on an Abelian surface,
and $\Gamma(\HMB)$ is $4$-dimensional.

If $s$ is a non-zero section of $\HMB$, and $p=[s]$ the corresponding point
in $\PS(\Gamma(\HMB))$, we will write $A_p$ for the zero locus of $s$.

\begin{definition}
Let $s_1$ and $s_2$ be two linearly independent sections of $\HMB$.
Then $s_1\wedge s_2$ is a non-zero section of
$\wedge^2\HMB\cong\OB(5)$,
and its zero-locus is a quintic threefold $V$;
we call such a threefold a \emph{\HM\ quintic}, or an \emph{HM quintic}.
\end{definition}

The threefold $V$ is determined by the line $L$
through the two points $[s_1]$ and $[s_2]$ in $\PS(\Gamma(\HMB))$,
and is covered by the pencil of surfaces
$\{\,A_p\mid p\in L\,\}$, giving a rational map $V\ratto L$.
We can resolve the indeterminacy of this rational map by taking
\begin{equation}\label{eq:defvtilde}
\~V=\{\,(x,p)\in V\times L\mid x\in A_p\,\}.
\end{equation}

In \cite[\secsign5]{HorMum73}, it is shown that the fivefold
\begin{equation*}
Z\defeq\{\,(x,p)\in\PS^4\times\PS(\Gamma(\HMB))\mid x\in A_p\,\}
\end{equation*}
is non-singular.
If $P_1$ and $P_2$ are two planes in $\PS(\Gamma(\HMB))$ such that
$L=P_1\cap P_2$, then $\~V=Z\cap(\PS^4\times P_1)\cap(\PS^4\times P_2)$,
and by Bertini's theorem, this is non-singular for generic
$P_1$ and $P_2$, so $\~V$ is non-singular for
generic $L$.

If $x\in V$, then the projection $\~V\to V$
is an isomorphism at
$x$ unless both $s_1$ and $s_2$ vanish at $x$, in which case
$\{x\}\times L\subset\~V$ is a rational curve contracted by the
projection.
Thus the singular locus of $V$ is $A_{[s_1]}\cap A_{[s_2]}$, which
for generic $s_1$ and $s_2$ consists of $100$ nodes.

\begin{definition}
Let $\epsilon=e^{2\pi i/5}$, and define three elements of $\SL_5\C$ by
\begin{equation*}
(\sigma x)_i=x_{i+1},\qquad
(\tau x)_i=\epsilon^ix_i,\qquad
(\iota x)_i=x_{-i}.
\end{equation*}
The \emph{Heisenberg group} $\HGp$ is the group of order $125$
generated by $\sigma$ and $\tau$.
\end{definition}

The centre of $\HGp$ is $\muGp$,
the group of $5$th roots of $1$, and $\HGp/\muGp\cong\Z_5^2$.
Recall from \cite{HorMum73} that the action of
$\langle\HGp,\iota\rangle$ on $\PS^4$ preserves each surface $A_p$,
and so preserves $V$ and can be lifted to an action on $\~V$.
Since $\muGp$ acts trivially on $\PS^4$, the $100$ nodes of $V$ form
$4$ $\HGp$-orbits of $25$ points, which pair up into
$2$ $\langle\HGp,\iota\rangle$-orbits of $50$ points.

\section{Determinantal quintics}\label{sec:detquintics}

This section considers determinantal quintics, of which
\HM\ quintics are a special case.

Let $(a_{ijk})$ be complex numbers, for $i,j,k\in\Z_5$,
and consider the matrix of linear forms
\begin{equation*}
M_{ij}(x)=\sum_ka_{ijk}x_k.
\end{equation*}
If $\det M(x)$ does not vanish identically,
the equation $\det M(x)=0$ defines
a quintic threefold $W$ in $\PS^4$.
Say that $M$ is \emph{nice} if $W$ is regular in codimension $1$,
$\rank M(x)\ge3$ everywhere, and
$\rank M(x)=4$ at all but finitely many points of $W$.

Let
\begin{equation*}
X=\{\,(x,y)\in\PS^4\times\PS^4\mid\forall i,\,
\sum_{j,k}a_{ijk}x_ky_j=0\,\}.
\end{equation*}
The projection $\pi_1$ takes $X$ to $W$,
with $\pi_1^{-1}\{x\}=\{x\}\times\PS(\ker M(x))$, so
if $M$ is nice, $\pi_1$ is a birational morphism,
contracting finitely many rational curves on $X$.

We can construct two more matrices of linear forms from $(a_{ijk})$:
\begin{gather*}
M'_{ij}(x)=\sum_ka_{jki}x_k,\\
M''_{ij}(x)=\sum_ka_{kij}x_k,
\end{gather*}
and the corresponding threefolds $W'$, $X'$, $W''$, and $X''$.
Note that $W'$ is the image of $X$ under the other projection
$\pi_2$; if $M'$ and $M''$ are also nice,
we get a sequence of birational maps,
all isomorphisms in codimension $1$:
\begin{equation}\label{eq:wxseq}
\newcommand{\sw}{\!\!{}^{\pi_1}\!\!\swarrow}
\newcommand{\se}{\searrow\!\!{}^{\pi_2}\!\!}
\newcommand*{\rt}[1]{\stackrel{#1}{\rattogets}}
\begin{array}{c*{12}{c@{}}cc}
&&&X&&&&X'&&&&X''&&&\\
\cdots&&\sw&&\se&&\sw&&\se&&\sw&&\se&&\cdots\\
&W&&\rt{\phi}&&W'&&\rt{\phi'}&&W''&&\rt{\phi''}&&W&
\end{array}
\end{equation}

Let $H$, $H'$, and $H''$ be hyperplane sections of $W$, $W'$,
and $W''$ respectively.
We can define an infinite sequence of divisors on $W$:
{
\newlength{\eqlen}\settowidth{\eqlen}{${}={}$}
\newcommand{\aldots}{\makebox[\eqlen]{$\vdots$}}
\begin{align*}
&\aldots\\
D_{-2}&=\phi''_*\phi'_*H'\\
D_{-1}&=\phi''_*H''\\
D_0&=H\\
D_1&=\phi^*H'\\
D_2&=\phi^*\phi^{\prime*}H''\\
&\aldots
\end{align*}
}

\begin{lemma}\label{lemma:deldelh}
For any $i\in\Z$, the divisor
$D_{i-1}+D_{i+1}$ is linearly equivalent to $4D_i$.
All the divisors $\pi_1^*D_i$ on $X$ are Cartier.
The cone in $\Cartsp(X)$ generated by $\pi_1^*D_i$ and $\pi_1^*D_{i+1}$
is the nef cone of a marked model of $X$,
isomorphic to $X$, $X'$, or $X''$ (according to the value of
$i$ modulo $3$).
\end{lemma}

\begin{proof}
(Here $\sim$ denotes linear equivalence.)
In \cite[p.~24]{Borcea91}, it is shown that $D_{-1}+D_1$
(there called $\Delta+\Delta'$) is the intersection of $W$ with a
quartic.  Therefore $D_{-1}+D_1\sim4H=4D_0$.
Applying the same result at other points in the sequence
(\ref{eq:wxseq}) gives $D_{i-1}+D_{i+1}\sim4D_i$ for all $i$.

Since $\pi_1^*D_0=\pi_1^*H$ and $\pi_1^*D_1=\pi_2^*H'$ are Cartier,
this equation then implies that $\pi_1^*D_i$ is Cartier for all $i$.

$\pi_1^*H$ and $\pi_2^*H'$
generate the nef cone of $X$; again, we can translate this
along the sequence (\ref{eq:wxseq}) to give the last statement of the
theorem.
\end{proof}

Solving the difference equation $D_{i-1}-4D_i+D_{i+1}=0$, we get
\begin{equation*}
D_i=\half\big(D_0+\tfrac{1}{\sqrt{3}}(D_1-2D_0)\big)(2+\sqrt{3})^i
+\half\big(D_0-\tfrac{1}{\sqrt{3}}(D_1-2D_0)\big)(2-\sqrt{3})^i,
\end{equation*}
so the limits of the rays generated by $D_i$ as $i\to\pm\infty$ are
generated by $D_0\pm\frac{1}{\sqrt{3}}(D_1-2D_0)$; these rays
then constitute the boundary of $\Mclone(X)\cap\Cartsp(X)$.
(For generic $M$, $X$ is smooth, so $\Cartsp(X)=\Weilsp(X)$ and
these rays are the boundary of $\Mclone(X)$.)

\medskip
Now let $y\in\C^5$, and consider the matrix $M_y$, where
\begin{equation*}
(M_y)_{i,j}(x)=x_{i+j}y_{i-j},\quad i,j\in\Z_5,
\end{equation*}
and the corresponding varieties and maps
\begin{equation}\label{eq:wyseq}
\cdots W_y\stackrel{\phi_y}{\rattogets}W'_y\stackrel{\phi'_y}{\rattogets}
W''_y\stackrel{\phi''_y}{\rattogets}W_y\cdots.
\end{equation}
Because $M_y^\transpose=M_{\iota y}$, this sequence can also be written as
\begin{equation}\label{eq:wiotayseq}
\cdots W_{\iota y}\stackrels{}{\phi''_{\iota y}}{\rattogets}
W''_{\iota y}\stackrels{}{\phi'_{\iota y}}{\rattogets}
W'_{\iota y}\stackrels{}{\phi_{\iota y}}{\rattogets}W_{\iota y}\cdots.
\end{equation}

From \cite[attributed there to R. Moore]{Borcea91} we see that,
if $V$ is a sufficiently general 
HM quintic, and $y$ is any node of $V$, then $V=W_y$; also
that there is an automorphism of the space of HM quintics, $\beta$,
such that $\beta V=W'_y$.

Therefore, (\ref{eq:wyseq}) and (\ref{eq:wiotayseq}) can also be written as
\begin{equation}\label{eq:vbetavseq}
\cdots V\stackrels{\phi_y}{\phi''_{\iota y}}{\rattogets}
\beta V\stackrels{\phi'_y}{\phi'_{\iota y}}{\rattogets}
\beta V\stackrels{\phi''_y}{\phi_{\iota y}}{\rattogets}V\cdots.
\end{equation}

\section{The small resolutions of $V$}\label{sec:smallresolutions}

From now on all threefolds we consider
will be isomorphic in codimension $1$
to a fixed sufficiently general HM quintic $V$.
We can therefore identify corresponding divisors
in any two of these threefolds.
Any such divisor determines a point in the real vector space
$N=\Weilsp(V)=\Cartsp(\~V)$.

Choose representatives $y_1,\ldots,y_4$ (respectively, $z_1,\ldots,z_4$) of
the four $\HGp$-orbits of nodes of $V$ (respectively, of $\beta^{-1}V$),
such that $\iota y_1=y_3$, $\iota y_2=y_4$, $\iota z_1=z_3$, and
$\iota z_2=z_4$.
Then $V$ fits into six different sequences of the form (\ref{eq:vbetavseq}):
\begin{gather*}
\cdots\beta V\stackrels{\phi''_{y_i}}{\phi_{\iota y_i}}{\rattogets}
V\stackrels{\phi_{y_i}}{\phi''_{\iota y_i}}{\rattogets}
\beta V\cdots,\quad i=1,2;\\
\cdots\beta^{-1}V\stackrels{\phi_{z_i}}{\phi''_{\iota z_i}}{\rattogets}
V\stackrels{\phi'_{z_i}}{\phi'_{\iota z_i}}{\rattogets}
V\cdots,\quad i=1,\ldots,4.
\end{gather*}
This gives us $12$ divisors on $V$:
\begin{itemize}
\item[]$\Delta_i$, the strict transform under $\phi_{y_i}$
of a hyperplane section of $\beta V$.
\item[]$\nabla_i$, the strict transform under
$\phi''_{\iota z_i}$ of a hyperplane section of $\beta^{-1}V$.
\item[]$\nabla'_i$, the strict transform under
$\phi_{z_i}$ of a hyperplane section of $V$.
\end{itemize}
We also have $H$, a hyperplane section of $V$, and $A$, one of the pencil
of Abelian surfaces on $V$.
By Lemma~\ref{lemma:deldelh}, we have
$\Delta_i+\iota\Delta_i\sim\nabla_i+\nabla'_i\sim4H$.

We write $\Lambda_i$ for the rational curve on $\~V$ that is contracted to
the node $y_i$ on $V$.  Noting from \cite{BHM87} that a singular member of
the pencil of surfaces $\{A_p\}$ is an elliptic scroll, we let $\Gamma$ be a
line in one of these scrolls.

The following intersection numbers on $\~V$ are computed in \cite{Borcea91}
($\Lambda_1, \Lambda_2, \Lambda_3, \Lambda_4$ are there called
$L_1, L_3, L_2, L_4$---note the different order):
\begin{equation}\label{eq:intersections}
\begin{array}{c*{14}{c}}
&H&A&\Delta_1&\Delta_2&\Delta_3&\Delta_4&
\nabla_1&\nabla_2&\nabla_3&\nabla_4&
\nabla'_1&\nabla'_2&\nabla'_3&\nabla'_4\\
\Gamma&    1&0& 2&2&2&2&  3&3&3&3&  1&1&1&1\\
\Lambda_1& 0&1& -1&0&1&0& 0&0&-1&-1& 0&0&1&1\\
\Lambda_2& 0&1& 0&-1&0&1& -1&0&0&-1& 1&0&0&1\\
\Lambda_3& 0&1& 1&0&-1&0& -1&-1&0&0& 1&1&0&0\\
\Lambda_4& 0&1& 0&1&0&-1& 0&-1&-1&0& 0&1&1&0
\end{array}\end{equation}

\begin{corollary}
$\{H,A,\Delta_1,\Delta_2\}$ is a basis for $N$;
the group $\HGp$ acts trivially on $N$; in
particular, the classes of $\Delta_i$ and
$\Lambda_i$ do not depend on the choice
of $y_i$ within its $\HGp$-orbit.
The class of $\Gamma$ does not depend on
which of the singular fibres it is chosen from.
\end{corollary}

\begin{proof}
In \cite[Proposition 2.2]{Aure89} it is proved that $\dim N=4$.
The multiplication table above shows that $H$, $A$,
$\Delta_1$, and $\Delta_2$ are linearly independent, so they
form a basis.  All are fixed by $\HGp$, so $\HGp$ acts trivially.
The intersections
of $\Gamma$ with these basis elements
do not depend on the singular fibre used,
so the numerical class of $\Gamma$ does not.
\end{proof}

\begin{convention}
Following \cite{Borcea91}, we will use the basis
\begin{equation}\label{eq:basis}
\{A,2H-\Delta_1,2H-\Delta_2,H\}
\end{equation}
for $N$.  The dual basis for $N^\dual$ is
$\{\half(\Lambda_1+\Lambda_3),\half(\Lambda_1-\Lambda_3),
\half(\Lambda_2-\Lambda_4),\Gamma\}$.
The components of a divisor $D$ will be written $(d_1,d_2,d_3,d_4)$.
\end{convention}

\begin{corollary}\label{cor:trivialautos}
The birational automorphisms of $V$ that act trivially on $N$ are exactly
the elements of $\HGp/\muGp$.
\end{corollary}

\begin{proof}
Any such automorphism fixes $H$, and so comes from
a projective automorphism of $V\subset\PS^4$.  These are exactly
$\langle\HGp,\iota\rangle/\muGp$, and $\iota$ does not act trivially.
\end{proof}

\begin{figure}[htp]
\includegraphics{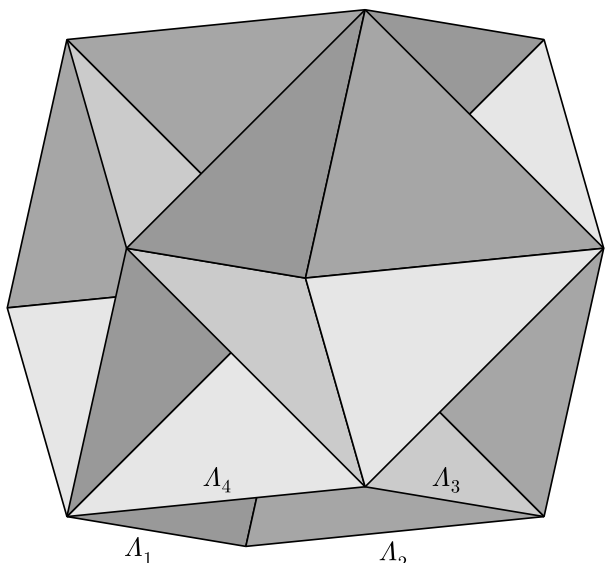}
\caption{\label{fig:cuboctos}
The vertex figure of an affine slice of the movable fan at $H$.}
\end{figure}

\begin{figure}[htp]
\includegraphics{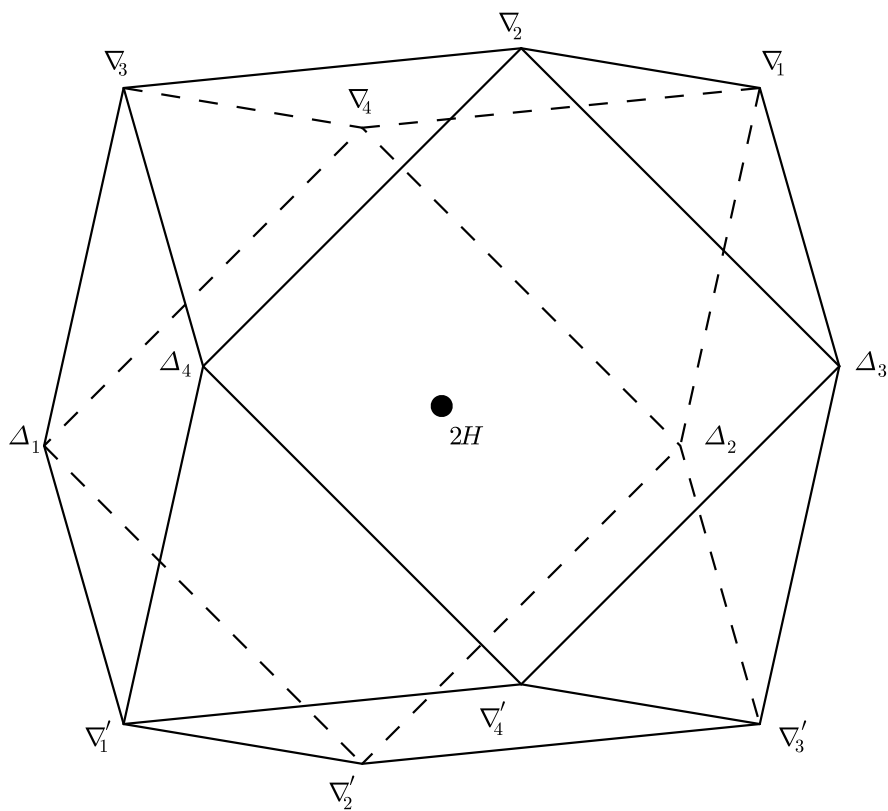}
\caption{\label{fig:cuboctof}
The nearest neighbours of $H$.}
\end{figure}

We can now draw a picture of a neighbourhood of the
ray generated by $H$ in $N$ (Figure~\ref{fig:cuboctos}),
showing the structure of the
movable fan there: taking the affine slice of the fan given by
$(\Gamma+\Lambda_1+\Lambda_3)\cdot D=2$, we find that
the four planes given by $\Lambda_i$ divide
the neigbourhood up in such a way that the vertex figure is a cuboctohedron
(as do any four planes in general position through a point).
The $12$ vertices of this cuboctohedron are on the lines joining $2H$
to the $12$ divisors constructed out of the determinantal description
($\Delta$s and $\nabla$s) (see Figure~\ref{fig:cuboctof}).

Let $\nabla$ be any of these $12$ divisors.  Coming as it does from
pulling back an ample divisor on a model of $V$ with Picard number $1$,
it generates a $1$-dimensional cone---a ray---in the movable fan.
(\ref{eq:wxseq}) in section~\ref{sec:detquintics} shows that
there is a marked model, not necessarily minimal,
of $V$, that is simultaneously a partial small resolution of
the two marked models corresponding to $H$ and $\nabla$, so there
is a cone in the fan whose closure contains both rays.

Therefore the $12$ rays generated by the $\Delta$s and $\nabla$s are
neighbours in the fan to the ray generated by $H$, and the
cuboctohedral arrangement extends as far as these neighbouring rays
(Figure~\ref{fig:cuboctof}).  Since
the marked models corresponding to these rays are all isomorphic,
ignoring the marking, either to $V$ or to $\beta V$,
there will be a similar cuboctohedral figure
about each of them, and so about every ray of the fan
that can be reached from
$H$ along $2$-dimensional cones of the fan.

In the next two sections we will describe
the nef cones of the resolutions of $V$---that is,
the $4$-dimensional
cones of the fan whose closure contains $H$.  This will
lead to a complete description of the movable fan.
These nef cones are the duals to the cones of curves described in
\cite{Borcea91}, but our notation and point of view are different,
and we need several details not found there.

Whenever we refer to the cuboctohedral neighbourhood of a ray in
the fan, we will use terms such as \emph{top} and \emph{left}, assuming
that the cuboctohedron is oriented as in Figures \ref{fig:cuboctos}
and~\ref{fig:cuboctof}.

\section{The Nef Cone of $\~V$.}\label{sec:tile}

The nef cone of $\~V$ lies in the sector of the cuboctohedral figure
given by the conditions $\Lambda_i\cdot D\ge0$; it contains
the square face of the cuboctohedron spanned by $\nabla'_1$,
$\nabla'_2$, $\nabla'_3$, and
$\nabla'_4$---that is, the bottom square face in the diagram.

In \cite[p.~31]{Borcea91} it is shown that $\phi'_{z_i}$ (there called
$\phi_{\nabla'_i}$) lifts to a biregular automorphism of $\~V$, acting by
translation on each Abelian surface $A_p$, and that $\phi'_{z_1}$ and
$\phi'_{z_2}$ generate a free abelian group of rank $2$; in terms of the basis
(\ref{eq:basis}), this group acts as the group of matrices
\begin{equation}
\left\{\,\begin{pmatrix}
1&2x&2y&x^2+y^2\\
0&1&0&x\\
0&0&1&y\\
0&0&0&1
\end{pmatrix}\,\Biggr|\,2x,2y,x-y\in\Z\,\right\}.
\end{equation}

Now we have enough information to describe the nef cone of $\~V$.

\begin{figure}[htp]
\includegraphics[scale=.75]{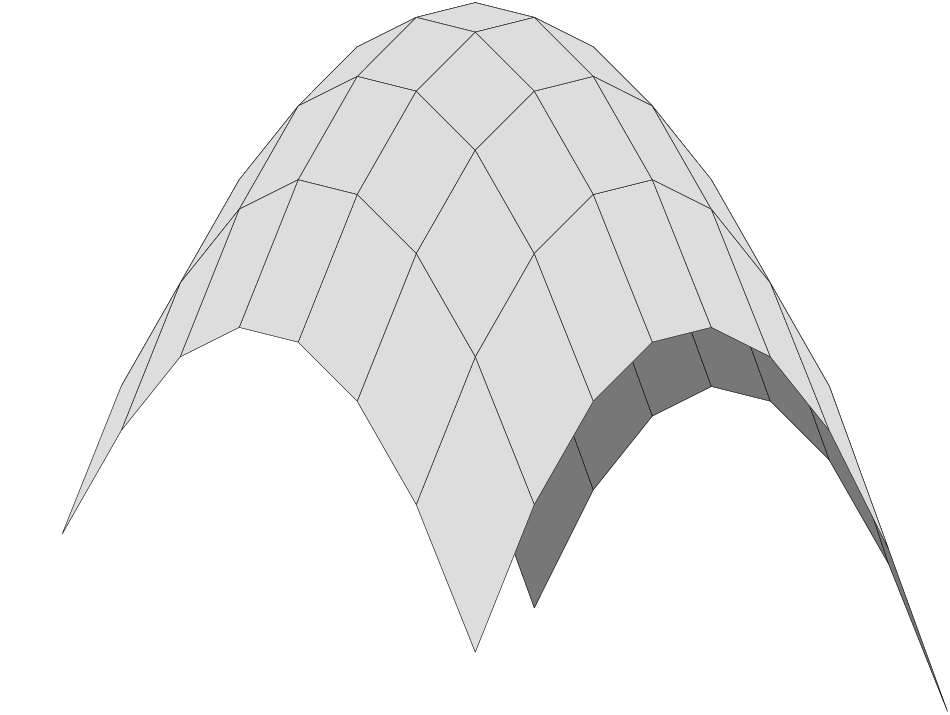}
\caption{\label{fig:tiling}
An affine slice of $\Aclone(\~V)$.}
\end{figure}

\begin{theorem}\label{thm:tile}
(i) $\Aclone(\~V)$ is generated by $A$ and all classes of the form
$(x^2+y^2,x,y,1)$ with $2x$, $2y$, and $x-y$ all integers.
In the affine slice given by $d_4=1$, this is a regular
square tiling of the paraboloid $d_1-d_2^2-d_3^2=0$ (see
Figure~\ref{fig:tiling}).

(ii) The automorphism group of $\~V$
acts transitively on the vertices and on
the faces of this tiling.

(iii) $(\Aut\~V)/(\HGp/\muGp)$ is isomorphic to the group
of all translation and order-$2$ rotation symmetries of this tiling.
\end{theorem}

\begin{proof}
(i) The set of vertices given is the
orbit of $H$ under the group generated by $\phi'_{z_1}$
and $\phi'_{z_2}$.
We described one of the faces (that given by $\Lambda_1$) above;
taking this information around
by the action of the group gives the whole tiling.
This tells us everything about the cone within the convex
hull of the tiling, and on the other hand the boundary of
the tiling is just the ray generated by $A$, which certainly
is nef, so
$\Aclone(\~V)$ is as claimed.

(ii) follows from the description of the cone in (i)---the translations
$\phi'_{z_i}$ are enough to give transitivity of the action on faces and
on vertices.

(iii) Any automorphism $\eta$
of $\~V$ must preserve the pencil of Abelian surfaces,
since it must, from the shape of the cone, preserve the class $A$.
It must preserve the individual surfaces in that pencil,
since the general member of
the pencil is not isomorphic to any other member.
Since the general surface in the pencil
has no automorphisms other than those
all Abelian surfaces have---translation
and reflection about a point---$\eta$ must act
in one of these two ways.
If the former, it commutes with other translations,
in particular $\psi_{z_1}$ and $\psi_{z_2}$, and so acts on the
tiling as a translation.  If the latter, $\eta\,\iota$ is a translation.
This shows that $(\Aut\~V)/\{\,$automorphisms acting trivially on
$N\,\}$ is generated by translations and $\iota$.
But by Corollary~\ref{cor:trivialautos}, the automorphisms that act trivially 
on $N$ are exactly $\HGp$.
\end{proof}

\begin{corollary}\label{cor:cubic}
The cubic intersection form on $\~V$ can be written as
\begin{equation}
D^3=5\Gamma\cdot D\,((\Gamma\cdot D)^2+3Q(D,D)),
\end{equation}
where $Q(D,E)=d_1e_4+d_4e_1-2d_2e_2-2d_3e_3$
in the basis~(\ref{eq:basis}).
\end{corollary}

\begin{proof}
The cubic form must be preserved by all automorphisms, and
$\Gamma$ and $Q$
generate the ring of
polynomial forms preserved by the automorphism group, so
$D^3$ must be a linear combination of $\Gamma^3$ and $\Gamma Q$.
Suppose $D^3=a\Gamma^3+b\Gamma Q$.

Then we can work out the values of $\Gamma$ and $Q$ on $H$ and
$A$, and
\begin{align*}
5=H^3
&=a(\Gamma\cdot H)^3+b(\Gamma\cdot H)Q(H,H)\\
&=a\ldot1^3+b\ldot1\ldot0=a,\\
10=H^2\cdot A
&=a(\Gamma\cdot H)(\Gamma\cdot A)^3
+\egypt b(\Gamma\cdot H)Q(A,H)
+\third b(\Gamma\cdot A)Q(H,H)\\
&=a\ldot1^2\ldot0+\egypt b\ldot1\ldot1+\third b\ldot0\ldot1=\egypt b,
\end{align*}
give $a=5$, $b=15$.
\end{proof}

\section{Second and Third Resolutions.}\label{sec:rest}

Since $\Aut\~V$ acts transitively on the faces of $\Aclone(\~V)$, all
models of $V$ that can be obtained by making one flop from $\~V$ are
isomorphic.
So let us consider the variety $\~V_1$ obtained by
flopping in the $\Lambda_1$ face,
so that $-\Lambda_1$ is effective on $\~V_1$.

\begin{figure}[htp]
\includegraphics[trim=0 10pt 0 10pt]{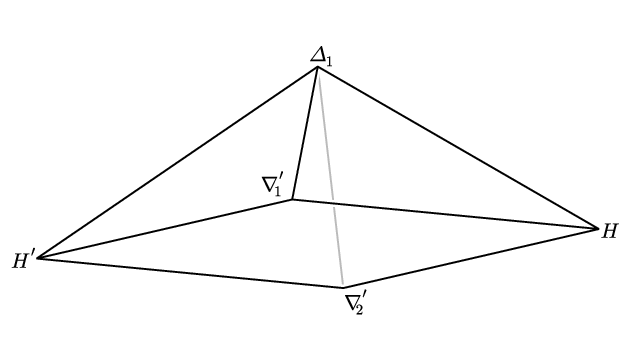}
\caption{\label{fig:pyramid}
An affine slice of $\Aclone(\~V_1)$.}
\end{figure}

The nef cone can be found easily in this case:

\begin{theorem}\label{thm:pyr}
An affine slice of $\Aclone(\~V_1)$
is the span of the five divisors
\begin{equation*}
H,\quad\nabla'_1,\quad
\nabla'_2,\quad H',\quad\Delta_1,
\end{equation*}
where $H'=\nabla'_1+\nabla'_2-H$.
It is a square pyramid, the base being the face at which
it meets $\Aclone(\~V)$ (see Figure~\ref{fig:pyramid}).
$(\Aut\~V_1)/(\HGp/\muGp)\cong\Z_2$, acting by a
rotation of order $2$ on the nef cone.
\end{theorem}

\begin{proof}
We are now considering the lower triangular faces of the cuboctohedron
of Figure~\ref{fig:cuboctof}.
Taking any one of the edges of the flopping face, the description of
the movable fan inside the cuboctohedron centred on that edge tells us that
half the boundary of the pyramid coincides with part of the boundary of
$\Aclone(\~V_1)$; these four descriptions put together
give the whole of the boundary of the pyramid, so the nef cone must
be the pyramid described.

The automorphism group must fix the only square face of the pyramid, so
$\Aut\~V_1$ consists exactly of those elements of $\Aut\~V$ that
fix the flopping face.
\end{proof}

The apex of the pyramid, $\Delta_1$, gives a contraction to $\beta V$,
and also has a cuboctohedral neighbourhood.  We can see that the cone
just described must fit into the top square face of this cuboctohedron, since
the bottom face is occupied by a model $\cong\beta\~V$, and the top
face is the only other face fixed by an automorphism of the movable fan of
order $2$.

Although there are two $\Aut\~V_1$-orbits of triangular faces of this
pyramid, the whole situation is symmetric (the two are swapped by changing
our arbitrary choice of $y_1$ and $y_2$ at the beginning of
section~\ref{sec:smallresolutions}),
so we only need to analyse one of the flops.
Let us consider $\~V_2$, obtained by flopping
$\~V$ in the planes defined by $\Lambda_1$ and $\Lambda_2$.

\begin{figure}[htp]
\includegraphics[trim=7pt 0 7pt 0]{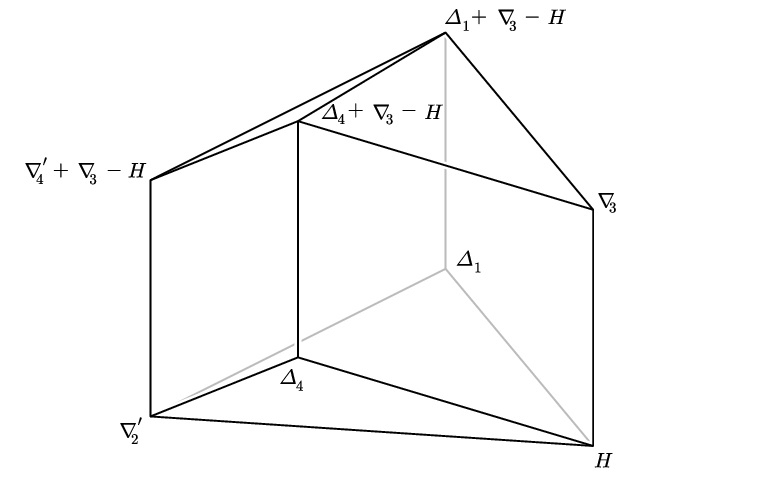}
\caption{\label{fig:lozengoid}
An affine slice of $\Aclone(\~V_2)$.}
\end{figure}

\begin{theorem}\label{thm:loz}
An affine slice of $\Aclone(\~V_2)$
is the span of the eight divisors
\begin{gather*}
H,\quad\Delta_1,\quad\Delta_4,\quad
\nabla'_1,\quad\nabla_3,\\
(\nabla'_1+\nabla_3-H),\quad
(\Delta_1+\nabla_3-H),\quad
(\Delta_4+\nabla_3-H).
\end{gather*}
We will call this shape a \emph{lozengoid}
(see Figure~\ref{fig:lozengoid}).
$(\Aut\~V_2)/(\HGp/\muGp)\cong\Z_2\times\Z_2$,
acting as the group of
orientation-preserving symmetries of this polyhedron.
\end{theorem}

\begin{proof}
We are now looking at the left-hand-side square face of the cuboctohedron
in Figure~\ref{fig:cuboctof}.  This immediately gives five
vertices:
$H$, $\Delta_1$, $\Delta_4$, $\nabla'_1$, and $\nabla_3$.
The cone also fits into a side square face of the equivalent
cuboctohedron centred at $\nabla'_1$, since $\~V_2$ is also
a resolution of $\nabla'_1$; there is a symmetry that interchanges
$H$ and $\nabla'_1$, which must also swap
$\Delta_1$ and $\Delta_4$, and this determines its effect
completely; it takes $\nabla_3$ to $\nabla'_2+\nabla_3-H$.

Now, looking from the point of view of the vertex $\Delta_1$,
the pyramidal cone of $\~V_1$, as discussed above, fits into the
top square face of the cuboctohedron, so $\~V_2$, one flop away,
fits into an upper triangular face.  Similarly, the pyramid that sits in
the top face of $H$'s cuboctohedron is in a lower triangular
face from the point of view of
$\nabla_3$, and $\~V_2$, two flops away, must again be in an
upper triangular face.  Symmetry shows that the same holds for
$\Delta_4$ and $\nabla'_2+\nabla_3-H$.
This means that we must have
a symmetry of the given slice of $\Aclone(\~V_2)$
any of these to any other,
which forces the cone to be the shape described.
\end{proof}

Since this involves the upper triangular faces of the cuboctohedron, it
completes the job of describing cones that meet $H$---that is,
nef cones of resolutions of $V$.
But now every face of any of these cones borders another cone
of one of these three types, 
so we have in fact described every cone in the fan:

\begin{theorem}
The minimal models of general HM quintics fall into three families,
the families of $\~V$, $\~V_1$, and $\~V_2$, that we will call,
after their nef cones, \emph{tiling type}, \emph{pyramid type} and
\emph{lozengoid type} respectively.
All the models satisfy Morrison's nef cone conjecture
\cite{Morrison93}; that is, the action of the biregular
automorphism group of each variety on its nef cone has
a rational polyhedral fundamental domain.

Any given HM quintic $V$ has infinitely
many marked minimal models in each family, but up to isomorphism has only
two minimal models of tiling type, two of pyramid type, and
four of lozengoid type.
\end{theorem}

\begin{proof}
Theorems \ref{thm:tile}, \ref{thm:pyr}, and~\ref{thm:loz} construct
the three types of cone, and the construction describes what is
on the other side of any face of any of these cones---there are
square faces where tilings meet pyramids, triangular faces where
pyramids meet lozengoids and quadrilateral faces where lozengoids
meet other lozengoids, twisted by a quarter-turn.
(This is
hard to visualize, since these faces are not square on the picture
(Figure~\ref{fig:lozengoid}), but of course these are \emph{affine} pictures
and we have chosen different affine slice for the different models.)
Every minimal model must be a finite sequence of flops away from
$\~V$, so we have described all minimal models.

The nef cone conjecture is clear, and a fundamental
domain is easily seen from the descriptions in the
three cases.

For the last statement, the two models of tiling type
are the resolutions of $V$ and $\beta V$, those
of pyramid type are one flop from each of these
(all such flops being isomorphic
because $\Aut(\~V)$ acts transitively on the faces of $\Aclone(\~V)$);
lozengoid-type models come from flopping
two adjacent faces of $\Aclone(\~V)$,
which can be chosen in two ways because of the lack of automorphisms
that act as order-$4$ rotations on the nef cone of $\~V$.
\end{proof}

\section{Combinatorics of the Movable Fan.}\label{sec:comb}

From the theorems of the last two sections we can deduce a complete
combinatorial
description of the movable fan of $V$, since they give the shapes
of the cones and the way the cones fit together at each ray.
In this section we aim to give an understanding of the combinatorics of
the movable fan as a polyhedral complex,
disregarding the linear structure, and of the action on that
complex of the birational automorphism group.
So we will allow ourselves to squash and fold the cones.
We will talk in terms of
affine slices, making it a problem about $3$-dimensional solids.

\begin{figure}[htp]
\includegraphics{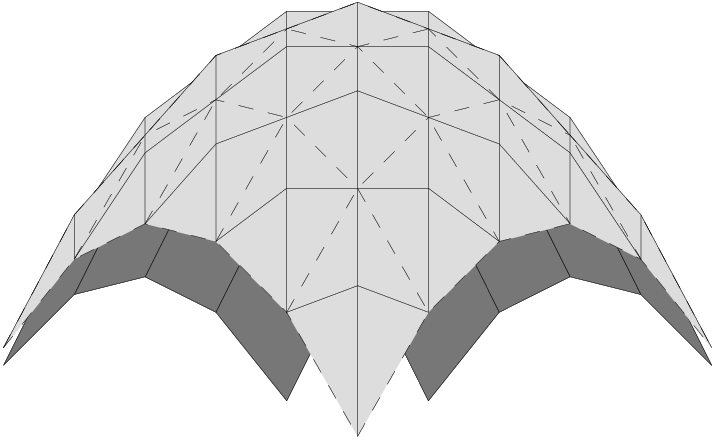}
\caption{\label{fig:tandp}
A stellation of a square tiling.}
\end{figure}

The first thing we can do to simplify matters is to glue the pyramids
to the faces of the tilings.
The result is a solid, no longer convex, with
twice the vertices of the old square tiling, arranged like a square
tiling again (at $45^\circ$ to the old one), but with half the
vertices protruding and with a diagonal inward fold
(marked by dotted lines in Figure~\ref{fig:tandp}) across each face.

\begin{figure}[htp]
\includegraphics[scale=.6]{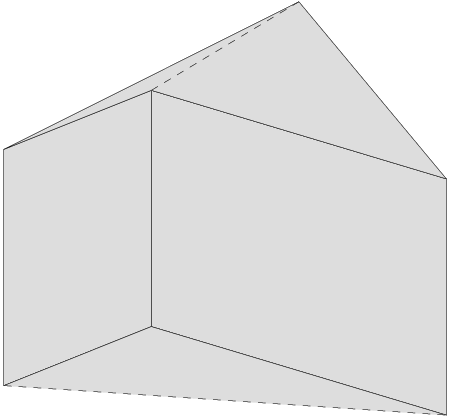}
\caption{\label{fig:lozascube} 
A lozengoid.}
\end{figure}

This reduces the problem to two sorts of solid, these stellated
tilings and the lozengoids.
Now the lozengoid of Figure~\ref{fig:lozengoid}
can be regarded as a cube distorted in such a
way that two opposite faces (the top and bottom faces 
in Figure~\ref{fig:lozengoid}) become folded outwards 
along their diagonals (see Figure~\ref{fig:lozascube}).

The way these cones are put together sets the dotted ridges on the
lozengoids into the dotted furrows on the tilings, so
we can flatten these folds to obtain a complex of ordinary cubes
and square-tiled paraboloids.

\begin{figure}[htp]
\includegraphics[scale=.5,trim=0 0 0 -12pt]{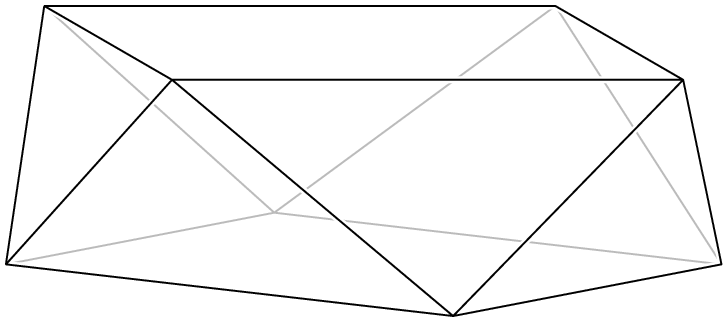}
\caption{\label{fig:antiprism}
A square antiprism.}
\end{figure}

Consider what has happened to the
vertex figure (Figure~\ref{fig:cuboctof}) in
this process:
gluing the pyramids to the tilings has the effect of
deleting the four edges around the bottom face of the cuboctohedron,
leaving four downward-pointing spikes, which are flattened out
when the dotted edges on the solids are flattened.
The result is a square antiprism (Figure~\ref{fig:antiprism}).
A square tiling pokes into the each of the two square faces
and a cube pokes into each triangular face.

It is easily seen that this vertex configuration determines a unique
simply connected polyhedral complex, because there is no choice
about how the figures at two adjacent vertices are
oriented with respect to one another.
But we can construct such a simply-connected polyhedral complex
with this configuration at every vertex as a \emph{uniform hyperbolic
tiling} (see \cite{Coxeter73}, \cite{Humphreys90}) as follows:

The Coxeter graph $\Cox$ describes a hyperbolic reflection group $C$
with a tetrahedral fundamental cell: the nodes of the graph
correspond to the faces of the cell, and the angle between two faces is
$\pi/8$ if the two nodes are joined by a double line, and $\pi/4$ if the
two nodes are not joined.
The marked graph $\MCox$ represents a hyperbolic tiling with
$C$ acting transitively on its vertices:
the unique vertex of the tiling that is in the fundamental cell
lies on the faces of the cell corresponding to unfilled nodes in
the marked graph, and not on the faces corresponding to filled nodes.
$C$-orbits of $k$-dimensional faces of the tiling are described 
by $k$-node subgraphs of $\MCox$ having at least one marked node 
in each connected component
(for example, deleting an end node gives $\TilCox$, a square tiling, 
and deleting a middle node gives $\CubCox$, a square prism), 
so it is possible to work out the vertex configuration of this
tiling: it is the same square antiprism as produced above.

Therefore we can identify 
the birational automorphism group of $V$:

\begin{theorem}
$\Bir(V)/(\HGp/\muGp)$ is isomorphic to a subgroup of index $8$ in 
the Coxeter group whose graph is $\Cox$.
\end{theorem}

\begin{proof}
Recall from Corollary~\ref{cor:trivialautos} that
$\HGp/\muGp$ is the group of birational automorphisms that act
trivially on $N$: therefore $\Bir(V)/(\HGp/\muGp)$ is
the group of all symmetries of the movable fan arising from
birational automorphisms.
The construction above identifies the group of all symmetries
of the fan with a subgroup of finite index in the group of symmetries
of the uniform tiling $T$ given by $\Cox$, which is larger than
the Coxeter group $C$ because the graph has a non-trivial automorphism.
But this graph automorphism exchanges the two $C$-orbits
of square tiling cells in $T$, and the construction above identifies
these two orbits with the class of marked models isomorphic to
$\~V$ and the class  of marked models isomorphic to $\beta\~V$,
so birational automorphisms do not swap them over.
Therefore $\Bir(V)/(\HGp/\muGp)$ is isomorphic to a subgroup of $C$.

In fact, the subgroup in question is of index $8$ in $C$---a factor of $2$
comes from the fact that the creased face diagonals
must be preserved (alternatively, because the two classes of
vertices corresponding to $V$ and $\beta V$ must be preserved)
and a further factor of $4$ because there are no
automorphisms of $\~V$ that induce
rotations of order $4$ or reflections.
\end{proof}

\section{The Boundary of the Movable Cone.}\label{sec:boundary}

The aim of this section is to demonstrate that, although by the
results of the last section the movable fan of $V$ is
combinatorially similar to a tiling of hyperbolic $3$-space,
the cone $\Mcone(V)$ itself is more subtle, by proving this
theorem:

\begin{theorem}\label{thm:notc2}
The boundary of the cone $\Mcone(V)$ is not twice differentiable
at any point in $\partial\Mcone(V)$.
\end{theorem}

The description of the nef cones in sections \ref{sec:tile} and~\ref{sec:rest}
allows us to find all the effective divisors on the boundary of the
movable cone, since every such divisor is nef on some model.
But the only rays in any of the nef cones described
above that are in the boundary of $\Mcone(V)$ are those corresponding
to Abelian fibrations on tiling type models.
These rays must then be dense on the boundary, and we can investigate the
shape of the boundary through them.  So let us see how we
can produce examples of these points.

From $\~V$ we can flop the four classes of curves
$\Lambda_1,\ldots,\Lambda_4$
to get a pyramid type model (the top of the cuboctohedron),
from which we can flop one more class, say $\Lambda'$, to get
a model $\~V'\cong\beta\~V$.  
We know that $\Lambda'\cdot\nabla_i=0$ for $i=1,\ldots,2$,
and $\Lambda'\cdot H=1$, which determine $\Lambda'$;
it can easily be checked that
$\Lambda'=3(\Lambda_1+\Lambda_3)+\Gamma$.

Let $A'$ be the class of an Abelian fibre
on $\~V'$, and $\Gamma'$ be the class of a line in a singular
fibre of the Abelian fibration given by $A'$.
Recall that in Corollary~\ref{cor:cubic} we defined an $\Aut\~V$-invariant
quadratic form $Q$ by $Q(D,E)=d_1e_4+d_4e_1-2d_2e_2-2d_3e_3$, in the
basis~(\ref{eq:basis}).  Let $Q'$ be the $\Aut\~V'$-invariant
quadratic form defined in the same way with reference to $\~V'$.

To get from $\~V$ to this model
we have the flopped the curves in each of the five classes
$\Lambda_1,\ldots,\Lambda_4$, and $\Lambda'$, so we know
$A'\cdot\Lambda_1=\cdots=A'\cdot\Lambda_4=A'\cdot\Lambda'=-1$,
which determines $A'$ exactly:
\begin{equation}
A'=5H-A=(-1,0,0,5).
\end{equation}

We can also do this process in reverse, flopping a face of $\Aclone(\~V)$,
say that given by $L_1$,
and four more faces to get a model $\~V''\cong\beta\~V$.
Then $\~V$ bears the relation to $\~V''$ that $\~V'$ bears to
$\~V$, so $A=5\Delta_1-A''$; that is,
\begin{equation}
A''=5\Delta_1-A=(-1,-5,0,10).
\end{equation}

Therefore we have
\begin{gather*}
Q(A',A')=-10,\qquad Q(A'',A'')=-70,\\
\Gamma\cdot A'=5,\qquad\Gamma\cdot A''=10,
\end{gather*}
so
\begin{align}
\label{eq:a'}
5Q(A',A')+2(\Gamma\cdot A')^2&=0,\\
\label{eq:a''}
5Q(A'',A'')+11(\Gamma\cdot A'')^2&=0.
\end{align}

This is enough to prove Theorem~\ref{thm:notc2}:

\begin{proof}
These two points $A'$ and $A''$
can be moved by the automorphisms of $\~V$ to get 
infinitely many points on the boundary of $\Mcone(V)$, accumulating
at $A$.
The equations (\ref{eq:a'}) and (\ref{eq:a''}) are preserved by the
automorphisms of $\~V$, so there are points accumulating at $A$
on two \emph{different} quadric cones tangent at $A$.

These two sets of points constrain the second
derivative, if it were defined, to take two different values,
so it cannot be defined.
\end{proof}

\begin{figure}[htp]
\includegraphics[scale=.7]{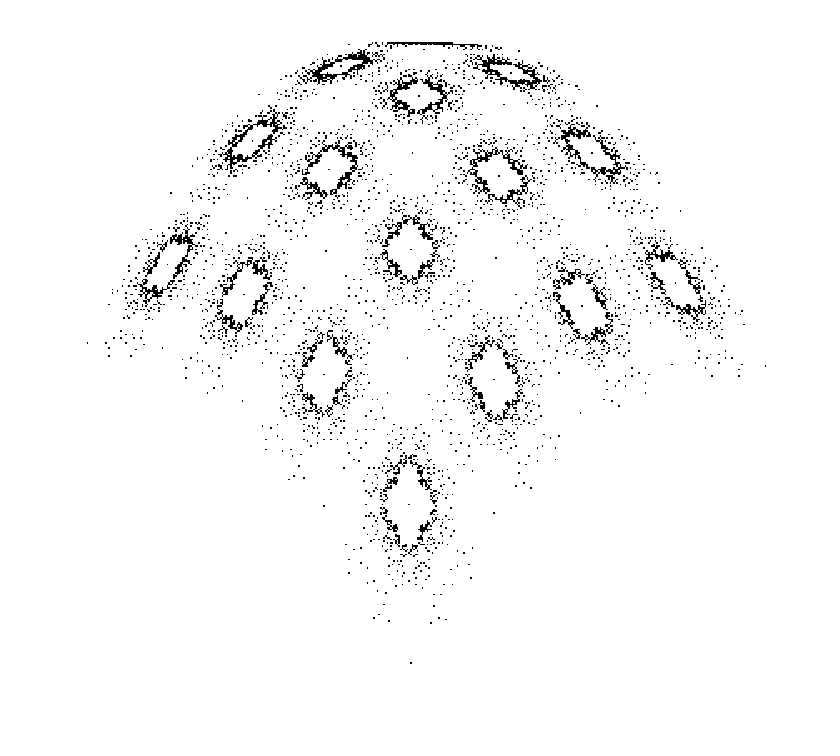}
\includegraphics[scale=.7]{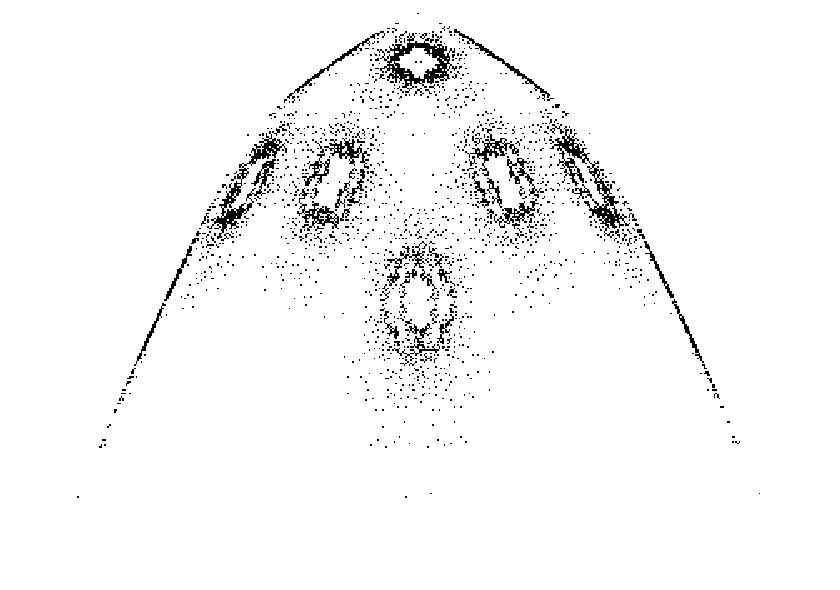}
\caption{\label{fig:dotplot}
Two views of part of the boundary of $\Mcone(V)$.}
\end{figure}

The two processes for obtaining new boundary rays from old can be iterated;
the plots in Figure~\ref{fig:dotplot} are of
two different orthogonal projections of and affine slice
of some such rays, showing the subtle shape of the boundary.

\section{The Movable Cone Conjecture.}

Finally, we can check explicitly that $V$ satisfies
Morrison's movable cone conjecture \cite{Morrison66}:

\begin{theorem}
There is a rational polyhedral fundamental domain for the action
of $\Bir(V)$ on $\Mcone(V)$.
\end{theorem}

\begin{proof}
First consider the hyperbolic uniform tiling constructed
in section~\ref{sec:comb} by distorting
the movable fan of $V$.
The standard fundamental domain for the action of the reflection group
on this tiling is a tetrahedron with a vertex
at the centre of one cell in each orbit---the centres of two cubes
and two square-tiled paraboloids
(where the `centre' of a paraboloid is the unique point
it contains that is in the boundary of hyperbolic $3$-space).
The thick-line triangles on Figure~\ref{fig:funddom}
show how eight such fundamental domains can be combined into a single
fundamental domain for the action of $\Bir(V)$.  

\begin{figure}[htp]
\includegraphics[scale=0.85]{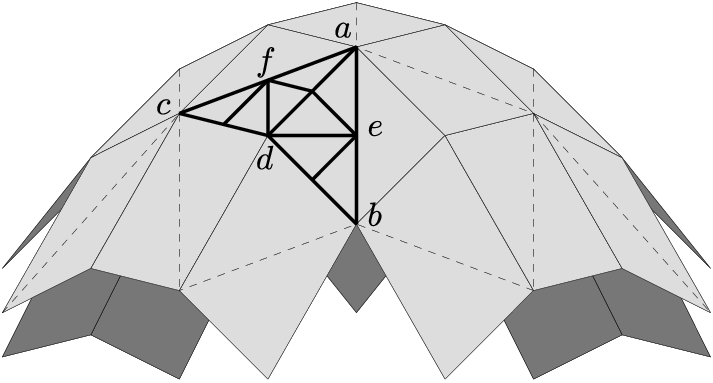}
\caption{\label{fig:funddom}
A fundamental domain for the action of $\Bir(V)$ on $\Mcone(V)$.}
\end{figure}

Its vertices are seven in number:
\begin{itemize}
\item[]$v_0$, the centre of the paraboloid drawn (in this view, the point
at infinity vertically downwards),
\item[]$v_a$, $v_b$, $v_c$, and $v_d$, the centres of the square-tiled
paraboloids meeting the one drawn at the points $a$, $b$, $c$, and $d$,
\item[]$v_e$ and $v_f$, the centres of the cubes meeting the paraboloid
drawn at $e$ and $f$.
\end{itemize}

The facets of this fundamental domain are the triangles
$v_bv_dv_e$, $v_av_dv_e$, $v_av_dv_f$, and $v_cv_dv_f$, and
the quadrilaterals $v_0v_bv_dv_c$, $v_0v_av_ev_b$, and $v_0v_cv_fv_a$.

Now consider the movable fan of $V$.  For each of the
points $v_i$, the centre of a square-tiled paraboloid or a cube, let
$\^v_i$ be a generator of the central ray of the corresponding
square-tiled paraboloid or lozengoid in the movable fan, all chosen
in the same affine slice $\Pi$.
It is clear from symmetry that the three quadrilaterals $\^v_0\^v_b\^v_d\^v_c$,
$\^v_0\^v_a\^v_e\^v_b$, and $\^v_0\^v_c\^v_f\^v_a$ are planar, and
the triangles $\^v_b\^v_d\^v_e$, $\^v_a\^v_d\^v_e$,
$\^v_a\^v_d\^v_f$, and $\^v_c\^v_d\^v_f$ are certainly planar,
so there is a (possibly non-convex) polyhedron with these seven faces.

The cone on this polyhedron will be a rational polyhedral fundamental domain
for the action of $\Bir(V)$ on $\Mcone(V)$ as long as it is convex.
Examination of the shape of the polyhedron shows that it will be convex as
long as the three quadrilaterals are convex and the line $\^v_e\^v_f$
passes through the interior of $\^v_0\^v_a\^v_d$.
We will therefore check these facts.

Choosing the right copy of the fundamental domain, we have, for some
positive scalars $\{\lambda_i\}$, the following expressions
(the expression for $\lambda_e\^v_e$
is the sum of the vertices of the lozengoid described in
Theorem~\ref{fig:lozengoid};
$A'$ was described in section~\ref{sec:boundary}, and
the remainder follow by translating these two by appropriate symmetries
of $\Aclone(\~V)$)
\begin{align*}
\lambda_0\^v_0&=A,\\
\lambda_a\^v_a&=A'=5H-A\\
\lambda_b\^v_b&=5\nabla'_1-A\\
\lambda_c\^v_c&=5\nabla'_2-A\\
\lambda_d\^v_d&=5\Delta_1-A
\end{align*}
\begin{align*}
\lambda_e\^v_e&=H+\Delta_1+\Delta_4+\nabla'_1+\nabla_3\\
&\qquad+(\nabla'_1+\nabla_3-H)+(\Delta_1+\nabla_3-H)+
(\Delta_4+\nabla_3-H)\\
&=10H+10\nabla'_1-6A\\
\lambda_f\^v_f&=10H+10\nabla'_2-6A
\end{align*}

Now the convexity of the quadrilaterals follows from
\begin{gather*}
\lambda_b\^v_b+\lambda_c\^v_c
=5\nabla'_1+5\nabla'_2-2A=5\Delta_1+3A=
\lambda_a\^v_a+4\lambda_0\^v_0,\cr
2\lambda_a\^v_a+2\lambda_b\^v_b=10H+10\nabla'_1-4A=
\lambda_e\^v_e+2\lambda_0\^v_0,\cr
2\lambda_a\^v_a+2\lambda_c\^v_c=10H+10\nabla'_2-4A=
\lambda_f\^v_f+2\lambda_0\^v_0,\cr
\end{gather*}
while the condition on $\^v_e\^v_f$ follows from
\begin{align*}
\lambda_e\^v_e+\lambda_f\^v_f&=
20H+10\nabla'_1+10\nabla'_2-12A\\
&=10\Delta_1+20H-2A\\
&=4\lambda_a\^v_a+2\lambda_d\^v_d+4\lambda_0\^v_0.
\end{align*}
All these equalities can be checked from the
table~(\ref{eq:intersections}).
\end{proof}

\clearpage
\providecommand{\bysame}{\leavevmode\hbox to3em{\hrulefill}\thinspace}

\bigskip
\begin{center}\small
Institut f\"ur Mathematik, Universit\"at Hannover,
Postfach~6009, D~30060 Hannover, Germany\\
$\langle\mbox{fryers@math.uni-hannover.de}\rangle$
\end{center}


\begin{thebibliography}{1}

\bibitem{Aure89}
A.~B. Aure, \emph{Surfaces on quintic threefolds associated to the
  {Horrocks--Mumford} bundle}, Arithmetic of Complex Manifolds, Lecture Notes
  in Math., no. 1399, Springer-Verlag, 1989, pp.~1--9.

\bibitem{BHM87}
W.~Barth, K.~Hulek, and R.~Moore, \emph{Degenerations of {Horrocks--Mumford}
  surfaces}, Math.\ Ann. \textbf{277} (1987), 735--755.

\bibitem{Borcea91}
C.~Borcea, \emph{On desingularized {Horrocks--Mumford} quintics}, J. Reine
  Angew.\ Math. \textbf{421} (1991), 23--41.

\bibitem{Coxeter73}
H.~S.~M. Coxeter, \emph{Regular polytopes}, third ed., Dover, 1973.

\bibitem{HorMum73}
G.~Horrocks and D.~Mumford, \emph{A rank 2 vector bundle on $\mathbb{P}^4$ with
  15,000 symmetries}, Topology \textbf{12} (1973), 63--81.

\bibitem{Humphreys90}
J.~Humphreys, \emph{Reflection groups and {Coxeter} graphs}, Cambridge Studies
  in Advanced Mathematics, no.~29, Cambridge University Press, 1990.

\bibitem{Kawamata88}
Y.~Kawamata, \emph{Crepant blowing-up of 3-dimensional canonical singularities
  and its application to degenerations of surfaces}, Ann.\ of Math.
  \textbf{127} (1988), 93--163.

\bibitem{Morrison66}
D.~Morrison, \emph{Beyond the {K\"ahler} cone}, Proceedings of the {Hirzebruch}
  65 Conference, Bar-Ilan Univ., 1966, pp.~361--376.
  alg-geom/9407007.

\bibitem{Morrison93}
\bysame, \emph{Compactifications of moduli spaces inspired by mirror symmetry},
  Journ\'ees de G\'eom\'etrie Alg\'ebrique d'{Orsay} 1992, Ast\'erisque, no.
  218, 1993, pp.~243--271.
  alg-geom/9304007.

\end{thebibliography}
\end{document}